\documentclass[12pt]{article}

\usepackage{amsfonts,amsmath,amssymb}
\usepackage{mathrsfs,mathtools,url}
\usepackage{latexsym}
\usepackage{tikz,color}

\newtheorem{theorem}{Theorem}[section]

\newtheorem{lemma}[theorem]{Lemma}
\newtheorem{proposition}[theorem]{Proposition}

\newtheorem{corollary}[theorem]{Corollary}

\let\deg\relax
\DeclareMathOperator {\deg} {deg}
\DeclareMathOperator {\BP} {BP}
\DeclareMathOperator {\bp} {bp}

\DeclareMathOperator {\muit} {\mu_{\rm it}}
\DeclareMathOperator {\mut} {\mu_{\rm t}}

\def\cp{\,\square\,}
\newcommand{\proof}{\noindent{\bf Proof.\ }}
\newcommand{\qed}{\hfill $\square$ \bigskip}

\vspace{8mm}
\title{Graphs with total mutual-visibility number zero and total mutual-visibility in Cartesian products}

\author{Jing Tian$\/^{a}$, Sandi Klav\v{z}ar$^{b, c, d}$  \\\\
$^{a}$ \small  College of Science, Nanjing University of
 Aeronautics \& Astronautics\\
 \small Nanjing, Jiangsu 210016, PR China\\
\small {\tt jingtian526@126.com}\\
$^{b}$ \small Faculty of Mathematics and Physics, University of Ljubljana, Slovenia\\
$^{c}$ \small Faculty of Natural Sciences and Mathematics, University of Maribor, Slovenia\\
$^{d}$ \small Institute of Mathematics, Physics and Mechanics, Ljubljana, Slovenia \\
\small {\tt sandi.klavzar@fmf.uni-lj.si}\\
}
\date{}

\begin{document}

\maketitle

\begin{abstract}
If $G$ is a graph and $X\subseteq V(G)$, then $X$ is a total mutual-visibility set if every pair of vertices $x$ and $y$ of $G$ admits a shortest $x,y$-path $P$ with $V(P) \cap X \subseteq \{x,y\}$. The cardinality of a largest total mutual-visibility set of $G$ is the total mutual-visibility number $\mu_{\rm t}(G)$ of $G$. Graphs with $\mu_{\rm t}(G) = 0$ are characterized as the graphs in which no vertex is the central vertex of a convex $P_3$. The total mutual-visibility number of Cartesian products is bounded and several exact results proved.  For instance, $\mu_{\rm t}(K_n\,\square\, K_m) = \max\{n,m\}$ and $\mu_{\rm t}(T\,\square\, H) = \mu_{\rm t}(T)\mu_{\rm t}(H)$, where $T$ is a tree and $H$ an arbitrary graph. It is also demonstrated that $\mu_{\rm t}(G\cp H)$ can be arbitrary larger than $\mu_{\rm t}(G)\mu_{\rm t}(H)$. 
\end{abstract}


\medskip\noindent
\textbf{AMS Math.\ Subj.\ Class.\ (2020)}: 05C12, 05C76


\section{Introduction}

Let $G = (V(G), E(G))$ be a graph and $X\subseteq V(G)$. Then $x, y\in X$ are {\em $X$-visible}, if there exists a shortest $x,y$-path $P$ with $V(P) \cap X = \{x,y\}$. $X$ is a \emph{mutual-visibility set} if its vertices are pairwise $X$-visible. The cardinality of a largest mutual-visibility set is the \emph{mutual-visibility number} $\mu(G)$ of $G$. A largest mutual-visibility set is a {\em $\mu$-set} of $G$.

These concepts were introduced and studied for the first time by Di Stefano in~\cite{distefano-2022}. The study was multiply motivated, notably by the role that mutual visibility plays in problems arising in the context of distributed computing by mobile entities, and by the fact that vertices in mutual visibility may represent entities on some nodes of a computer/social network that want to efficiently communicate in such a way that the exchanged messages do not pass through other entities. We should also mention related (but not identical) concepts in computer science that have been explored:  distributed computing by mobile entities~\cite{flocchini-2012}, the so-called mutual visibility task~\cite{diluna-2017}, and fat entities modelled as disks in the Euclidean plane~\cite{poudel-2019}. A related graph theory topic is also the general position in graphs, introduced in~\cite{manuel-2018, ullas-2016} and extensively studied by now, cf.~\cite{patkos-2020, yao-2022}. The problem was especially well investigated on Cartesian product graphs~\cite{ghorbani-2021, klavzar-2021a, klavzar-2021b, tian-2021a, tian-2021b}.

In~\cite{cicerone-2023}, the mutual-visibility problem was studied on Cartesian products and on triangle-free graphs, while in~\cite{cicerone-2023+} the focus was on strong products. In these studies, the following tools have proven to be extremely useful. We say that $X\subseteq V(G)$ is a {\em total mutual-visibility set} of $G$ if every pair of vertices $x$ and $y$ of $G$ is $X$-visible, that is, there exists a shortest $x,y$-path $P$ with $V(P) \cap X \subseteq \{x,y\}$. The cardinality of a largest total mutual-visibility set of $G$ is the {\em total mutual-visibility number} $\mut(G)$ of $G$. Further, $X$ is a {\em $\mut$-set} if it is a total mutual-visibility set of largest possible cardinality.

As observed in~\cite{cicerone-2023}, there exist graphs $G$ with $\mut(G) = 0$, that is, with no total mutual-visibility sets at all. Moreover, partial results on such graphs were proved, in particular cactus graphs $G$ with $\mut(G) = 0$ characterized. In Section~\ref{sec:zero} we characterize general graphs $G$ for which $\mut(G) = 0$ holds as the graphs that contain no bypass vertices. We introduce the latter concept in Section~\ref{sec:prelim}, where we also give further definitions needed, recall some know results, and add a few additional  preparatory results. In Section~\ref{sec:cp} we prove bounds  on the total mutual-visibility number of Cartesian product graphs and demonstrate their sharpness by several exact results. For instance, $\mut(K_n\cp K_m) = \max\{n,m\}$ and $\mut(T\cp H) = \mut(T)\mut(H)$, where $T$ is a tree. In Section~\ref{sec:cpII} we continue by the investigation of Cartesian products by considering the estimate $\mut(G\cp H) \le \mut(G) \mut(H)$. It holds in many cases, but on the other hand we show that $\mut(G\cp H)$ can be arbitrary larger than $\mut(G)\mut(H)$. We conclude by suggesting some open problems and directions for further investigation.

\section{Preliminaries and bypass vertices}
\label{sec:prelim}

We first recall some definitions, for all other undefined graph theory concepts, see~\cite{west-2001}. All the graph in this paper are simple and connected, unless stated otherwise. The order of graph $G = (V(G), E(G))$ is denoted by $n(G)$, the minimum degree of $G$ is denoted by $\delta(G)$, and the subgraph of $G$ induced by $S\subseteq V(G)$ is denoted by $G[S]$. A subgraph $H$ of $G$ is \textit{convex} if for each vertices $x,y\in V(H)$, all shortest $x,y$-paths in $G$ lie completely in $H$. The \textit{girth} $g(G)$ of a graph $G$ with a cycle is the length of its shortest cycle. $S\subseteq V(G)$ is an \textit{independent set} if $G[S]$ is an edgeless graph. The cardinality of a largest independent set is the \textit{independence number} $\alpha(G)$ of $G$. If $\tau$ and $\tau'$ are two graph invariants, then we say that a graph $G$ is a $(\tau,\tau')$-graph if $\tau(G) = \tau'(G)$.

The Cartesian product $G\cp H$ of graphs $G$ and $H$ has the vertex
set $V(G\cp H) = V(G)\times V(H)$, vertices $(g, h$) and $(g', h')$ are adjacent if either $gg'\in E(G)$ and $h = h'$, or $g = g'$ and $hh'\in E(H)$.
Given a vertex $h\in V(H)$, the subgraph of $G\cp H$ induced by the set $\{(g,h): g\in V(G)\}$, is a {\em $G$-layer}
and is denoted by $G^h$. $H$-layers $^gH$ are defined analogously. Each $G$-layer and each $H$-layer is isomorphic to $G$ and $H$, respectively. Moreover, it is also well-known that each layer of a Cartesian product is its convex subgraph. We will use this fact later on many times, sometimes implicitly.

We next recall some known results. 

\begin{proposition} {\rm \cite[Corollary~4.3]{distefano-2022}}\label{prop:for trees}
Let $T$ be a tree and $L$ the set of its leaves. Then $L$ is a mutual-visibility set and $\mu(T)=|L|$.
\end{proposition}

\begin{proposition} {\rm \cite[Proposition~3.3]{cicerone-2023+}}\label{prop:for block graphs}
Block graphs (and hence trees and complete graphs) and graphs containing a universal vertex are all $(\mu,\mut)$-graphs.
\end{proposition}

\begin{proposition} {\rm \cite[Proposition~3.1]{cicerone-2023+}}\label{prop:for convex subgraph}
Let $G$ be a graph. If $V(G)=\bigcup_{i=1}^{k} V_i$, where $G[V_i]$ is a convex subgraph of $G$ and $\mut(G[V_i])=0$ for each $i\in [k]$, then $\mut(G)=0$.
\end{proposition}

Here and later on,  $[k]$ stands for $\{1,\ldots,k\}$. The following fact will be used several times later on.

\begin{proposition}
\label{prop:subsets-are-ok}
If $X$ is a total mutual-visibility set of a graph $G$ and $Y\subseteq X$, then $Y$ is also a total mutual-visibility set of $G$.
\end{proposition}

\proof
Let $X$ be a total mutual-visibility set of $G$ and let $x\in X$. Then it suffices to verify that $Y = X\setminus \{x\}$ is also a total mutual-visibility set of $G$. Let $u, v\in V(G)$. If $u, v\in Y$, then $u, v\in X$ and hence they are $X$-visible and thus also $Y$-visible. Assume next $u\in Y$ and $v\notin Y$. Then no matter whether $v=x$ or $v\in V(G)\setminus X$, the vertices $u$ and $v$ are $X$-visible and so also $Y$-visible. The last case is when $u,v\in V(G)\setminus Y$, but also in this case we have the same conclusion.
\qed

To conclude the preliminaries we introduce the following concept which appears essential in investigation of the total mutual-visibility concept. We say that a vertex $u$ of a graph $G$ is a {\em bypass vertex} if $u$ is not the middle vertex of a convex $P_3$ in $G$. Otherwise we say that $u$ is a {\em non-bypass vertex}. Let $\BP(G)$ denote the set of bypass vertices of $G$ and let $\bp(G) = |\BP(G)|$. For instance, if $n\ge 1$, then $\BP(K_n) = V(K_n)$ because there are no convex paths $P_3$ in a complete graph. Hence $\bp(K_n) = n$. Similarly, $\bp(K_{n,m}) = n + m$ for $n,m\ge 2$. Indeed, if $u,v,w$ induce a $P_3$ in $K_{n,m}$, then since $n, m \ge 2$, there exists a common neighbor $v'$ of $u$ and $w$, where $v'\ne v$, hence no $P_3$ in $K_{n,m}$ is convex. On the other hand, if $n\ge 5$, then $\bp(C_n) = 0$.

The basic fact on bypass vertices is the following.

\begin{lemma}
\label{lem:non-bypass-in-not-in}
If $u$ is a non-bypass vertex of a graph $G$ and $X$ is a total mutual-visibility set of $G$, then $u\notin X$.
\end{lemma}

\proof
Since $u$ is a non-bypass vertex of $G$, it is the central vertex of a convex $P_3$. If $x$ and $y$ are the neighbors of $u$ on this $P_3$ and $u$ would lie in $X$, then $x$ and $y$ would not be $X$-visible. Hence $u\notin X$.
\qed

Lemma~\ref{lem:non-bypass-in-not-in} implies that
\begin{equation}
\label{eq:bounded-by-bp}
\mut(G) \le \bp(G)\,.
\end{equation}

This bound is sharp. If $T$ is a tree with $n(T) \ge 3$, then using Proposition~\ref{prop:for trees} we get $\mut(T)=\bp(T)$. Similarly, $\bp(K_n) = n$. On the other hand, consider complete bipartite graphs $K_{n,m}$, $n,m\ge 3$. From~\cite[Corollary 3.6]{cicerone-2023+} and~\cite[Theorem 4.9]{distefano-2022} we know that $\mut(K_{n,m}) = \mu(K_{n,m}) = n + m - 2$, but $\bp(K_{n,m}) = n + m$. The graph $G$ from Fig.~\ref{fig:graph-G} is another sporadic example for which the bound~\eqref{eq:bounded-by-bp} is strict.  We have $\mut(G) = 1$ and  $\bp(G) = 2$, where $\BP(G) = \{g_6, g_7\}$.

\begin{figure}[ht!]
\begin{center}
\begin{tikzpicture}[scale=1.0,style=thick]
\tikzstyle{every node}=[draw=none,fill=none]
\def\vr{3pt} 

\begin{scope}[yshift = 0cm, xshift = 0cm]
 \node [below=0.5mm] at (0,0) {$g_2$};
 \node [below=0.5mm] at (1.5,0) {$g_1$};
 \node [above=0.5mm] at (0,1.5) {$g_3$};
 \node [above=0.5mm] at (1.5,1.5) {$g_4$};
 \node [below=0.5mm] at (2,0.75) {$g_5$};
 \node [below=0.5mm] at (3.5,0.75) {$g_6$};
 \node [above=0.5mm] at (5,0.75) {$g_7$};
  \node [below=0.5mm] at (6.5,0.75) {$g_8$};
 \node [below=0.5mm] at (7,0) {$g_9$};
 \node [below=0.5mm] at (8.5,0) {$g_{10}$};
 \node [above=0.5mm] at (8.5,1.5) {$g_{11}$};
 \node [above=0.5mm] at (7,1.5) {$g_{12}$};

\path (0,0) coordinate (x2);
\path (1.5,0) coordinate (x1);
\path (0,1.5) coordinate (x3);
\path (1.5,1.5) coordinate (x4);
\path (2,0.75) coordinate (x5);
\path (3.5,0.75) coordinate (x6);
\path (5,0.75) coordinate (x7);
\path (6.5,0.75) coordinate (x8);
\path (7,0) coordinate (x9);
\path (8.5,0) coordinate (x10);
\path (8.5,1.5) coordinate (x11);
\path (7,1.5) coordinate (x12);
\draw (x1) -- (x2) -- (x3)--(x4) -- (x5)-- (x6) -- (x7)-- (x8)-- (x9) -- (x10)-- (x11)-- (x12)--(x8);
\draw (x1) -- (x5);
\draw (x5) -- (x7);
\draw (x6) -- (x8);
\draw (x1) -- (x5);
\draw (x5) .. controls (3.5,1.25).. (x7);
\draw (x6) .. controls (5,0.25).. (x8);

\draw (x1)  [fill=white] circle (\vr);
\draw (x2)  [fill=white] circle (\vr);
\draw (x3)  [fill=white] circle (\vr);
\draw (x4)  [fill=white] circle (\vr);
\draw (x5)  [fill=white] circle (\vr);
\draw (x6)  [fill=white] circle (\vr);
\draw (x7)  [fill=white] circle (\vr);
\draw (x8)  [fill=white] circle (\vr);
\draw (x9)  [fill=white] circle (\vr);
\draw (x10)  [fill=white] circle (\vr);
\draw (x11)  [fill=white] circle (\vr);
\draw (x12)  [fill=white] circle (\vr);

\end{scope}
\end{tikzpicture}
\end{center}
\caption{Graph $G$. }
\label{fig:graph-G}
\end{figure}

\section{Graphs with $\mut = 0$}
\label{sec:zero}

In this section we characterize graphs with $\mut = 0$ and give several applications of the characterization. We begin by two lemmas, the second one being also of independent interest.

\begin{lemma}
\label{lem:mut=0 - iff - singletons}
Let $G$ be a graph. Then $\mut(G) = 0$ if and only if for every $x\in V(G)$, the set $\{x\}$ is not a total mutual-visibility set of $G$.
\end{lemma}

\proof
If $\mut(G) = 0$ and $x\in V(G)$, then $\{x\}$ is clearly not a total mutual-visibility set of $G$. To prove the other direction assume that $\mut(G) > 0$. Let $X$ be an arbitrary non-trivial total mutual-visibility set of $G$. If $x\in X$, then by Proposition~\ref{prop:subsets-are-ok}, $\{x\}$ is a total mutual-visibility set of $G$.
\qed

\begin{lemma}
\label{lem:bypass-vertex-is-good}
Let $G$ be a graph with $n(G)\ge 2$ and $u\in V(G)$. Then $\{u\}$ is a total mutual-visibility set of $G$ if and only if $u$ is a bypass vertex.
\end{lemma}

\proof
First assume that $u$ is a bypass vertex of $G$ and let $P$ be an arbitrary shortest $x,y$-path. If $P$ does not contain $u$, then $x$ and $y$ are of course $\{u\}$-visible. Assume hence that $u\in V(P)$. Let $x'$ and $y'$ be the neigbors of $u$ on $P$. (It is possible that $\{x,y\} \cap \{x',y'\} \ne \emptyset$.) By our assumption, $x'-u-y'$ is not a convex $P_3$. Since $P$ is a shortest path this implies that there exists a vertex $u'\ne u$ adjacent to $x'$ and $y'$. But then by replacing the subpath $x'-u-y'$ by $x'-u'-y'$ we obtain a shortest $x,y$-path which does not pass $u$. Hence $x$ and $y$ are $\{u\}$-visible and so $\{u\}$ is a total mutual-visibility set.

Conversely, assume that $\{u\}$ is a total mutual-visibility set. Assume further that $u$ is the central vertex of a $P_3$, where $x$ and $y$ are the neighbors of $u$. As $x$ and $y$ are $\{u\}$-visible, we have an edge $xy$ or there is a vertex $u'$ adjacent to $x$ and $y$. Hence, the $P_3$ is not convex and so $u$ is a bypass vertex.  
\qed

Note that Lemma~\ref{lem:bypass-vertex-is-good} does not extend to two vertices. For instance, two opposite vertices of $C_4$ are bypass vertices, but they do not form a total mutual-visibility set.

The announced characterization now reads as follows.

\begin{theorem}
\label{thm:main-characterization-for-0}
Let $G$ be a graph with $n(G)\ge 2$. Then $\mut(G) = 0$ if and only if $\bp(G) = 0$.
\end{theorem}

\proof
If $\mut(G) = 0$, then $\bp(G) = 0$ by Lemma~\ref{lem:bypass-vertex-is-good}. Conversely, assume that $\bp(G) = 0$. By Lemma~\ref{lem:mut=0 - iff - singletons} it suffices to show that each set $\{x\}$ is not a total mutual visibility set of $G$. Suppose on the contrary that there is a vertex $x\in V(G)$ such that $X = \{x\}$ is a total mutual visibility set of $G$. By our assumption, $x$ is a non-bypass vertex, hence there exists a convex $P_3$ with $x$ the central vertex of it. But then the neihbors of $x$ on this $P_3$ are not $X$-visible, hence $X$ is not a total mutual-visibility set.
\qed

Clearly, to check whether a vertex is a bypass vertex is algorithmically simple. Hence the characterization of graphs $G$ with $\mut(G) = 0$ from Theorem~\ref{thm:main-characterization-for-0} is efficient.

Note that Theorem~\ref{thm:main-characterization-for-0} implies that if $\mut(G) = 0$, then $\delta(G) \ge 2$. Another consequence of the theorem is the following.

\begin{corollary}
\label{cor:girth}
Let $G$ be a graph with $g(G)\ge 5$. Then $\mut(G) = 0$ if and only if $\delta(G) \ge 2$.
\end{corollary}

\begin{figure}[ht!]
\begin{center}
\begin{tikzpicture}[scale=1.0,style=thick]
\tikzstyle{every node}=[draw=none,fill=none]
\def\vr{3pt} 

\begin{scope}[yshift = 0cm, xshift = 0cm]
\path (1,0) coordinate (x1);
\path (2.5,0) coordinate (x2);
\path (4,0) coordinate (x3);
\path (5.5,0) coordinate (x4);
\path (7,0) coordinate (x5);
\path (2.5,-1.5) coordinate (x6);
\path (5.5,-1.5) coordinate (x7);
\path (4,1) coordinate (x8);
\path (2.5,2) coordinate (x9);
\path (5.5,2) coordinate (x10);
\draw (x1) -- (x2) -- (x3) -- (x4) -- (x5);
\draw (x1) -- (x6) -- (x2);
\draw (x4) -- (x7) -- (x5);
\draw (x3) -- (x8) -- (x9);
\draw (x8) -- (x10);
\draw (x1) -- (x9) -- (x10) -- (x5) -- (x7) -- (x6) -- (x1);

\draw (x1)  [fill=white] circle (\vr);
\draw (x2)  [fill=white] circle (\vr);
\draw (x3)  [fill=white] circle (\vr);
\draw (x4)  [fill=white] circle (\vr);
\draw (x5)  [fill=white] circle (\vr);
\draw (x6)  [fill=white] circle (\vr);
\draw (x7)  [fill=white] circle (\vr);
\draw (x8)  [fill=white] circle (\vr);
\draw (x9)  [fill=white] circle (\vr);
\draw (x10)  [fill=white] circle (\vr);
\end{scope}
\end{tikzpicture}
\end{center}
\caption{A graph with the total mutual visibility number $0$.}
\label{fig:Halin-example}
\end{figure}
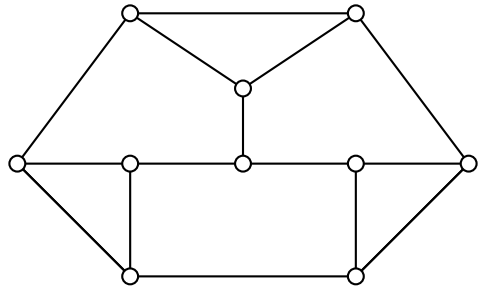

The Petersen graph applies to Corollary~\ref{cor:girth}. In addition, the corollary implies the characterization of cactus graphs $G$ with $\mut(G) = 0$ as given in~\cite[Proposition~3.2]{cicerone-2023+}. For a sporadic example of a graph $G$ with $\mut(G) = 0$ see Fig.~\ref{fig:Halin-example}.

As another application of Theorem~\ref{thm:main-characterization-for-0} we next determine the theta graphs with $\mut = 0$. For any positive integer $k\geq 2$ and $1\le p_1\leq\cdots \leq p_k$, the \textit{theta graph} $\Theta(p_1,\ldots, p_k)$ is the graph consisting of two vertices $a$ and $b$ which are joined by $k$ internally disjoint paths of respective lengths $p_1,\ldots, p_k$, where $p_2\ge 2$. (We add that several authors use the name theta graph restricted to the case $k=3$ in our definition, cf.~\cite{zhai-2021}.)

\begin{corollary}
\label{cor:theta}
If $1\le p_1\leq\cdots \leq p_k$, where $k\ge 2$ and $p_2\ge 2$, then $\mut(\Theta(p_1, \ldots, p_k)) = 0$ if and only if the following cases hold:
(i) $p_1=1$ and $p_2\geq 4$;
(ii) $p_1=2$ and $p_2\geq 3$;
(iii) $p_1\geq 3$.
\end{corollary}

\proof
Set $\Theta = \Theta(p_1, \ldots, p_k)$ and let $P_1, \ldots, P_k$ be the respective paths of $\Theta$ connecting $a$ and $b$. If $p_1 = 1$ and $p_2\le 4$, then an inner vertex of $P_2$ does not fulfill  the condition of Theorem~\ref{thm:main-characterization-for-0}, hence $\mut(\Theta) \ge 1$. On the other hand, if (i) $p_1=1$ and $p_2\geq 4$, then Theorem~\ref{thm:main-characterization-for-0} implies $\mut(\Theta) = 0$. The remaining cases are treated similarly.
\qed

We conclude the section with a description of Cartesian products with $\mut = 0$.

\begin{theorem}
\label{thm:cp}
If $G$ and $H$ are graphs, then $\mut(G\cp H) = 0$ if and only if $\mut(G) = 0$ or $\mut(H) = 0$.
\end{theorem}

\proof
Assume first that $\mut(G\cp H) = 0$. Suppose on the contrary that $\mut(G)\geq 1$ with a total mutual-visibility set $X$ and $\mut(H)\geq 1$ with a total mutual-visibility set $Y$.
By Proposition \ref{prop:subsets-are-ok},
there exist two vertices $x\in X$ and $y\in Y$ such that the sets $\{x\}$ and $\{y\}$ are total mutual-visibility set of $G$ and $H$.
Then we claim that $U=\{u\}$ with $u=(x,y)$ is a total mutual-visibility set of $G\cp H$.

Let $v, w$ be arbitrary vertices of $G\cp H$. We need to show that they are $U$-visible. If $v = u$ or $w = u$, there is nothing to be proved. If $v$ and $w$ lie in the $G$-layer containing $u$, then $v$ and $w$ are $U$-visible because their projections onto $G$ are $\{x\}$-visible. By the same argument we see that $v$ and $w$ are $U$-visible when $v$ and $w$ lie in the $H$-layer containing $u$. Assume hence that $v$ and $w$ neither lie in a common $G$-layer or a common $H$-layer. Then $w=(g_1,h_1)$ and $v=(g_2,h_2)$, where $g_1\ne g_2$ and $h_1\ne h_2$. Then it is well known that there exist two internally disjoint $v,w$-shortest paths. Since at least one of these two paths does not contain $u$, $v$ and $w$ are $U$-visible in $G\cp H$ also in this case. Hence we conclude that $\mut(G\cp H)\geq 1$.

To prove the converse, we may assume, without loss of generality, that $\mut(G)=0$.
Since $G^{h}$ is a convex subgraph of $G\cp H$ for any $h\in V(H)$, by Proposition~\ref{prop:for convex subgraph}, we have $\mut(G\cp H)=0$.
By symmetry, the same result also holds if $\mut(H)=0$, as desired.
\qed

Theorem \ref{thm:cp} extends to an arbitrary number of factors as follows.

\begin{corollary}
\label{cor:cp}
If $G = G_1\cp \cdots \cp G_k$, where $k\geq 2$, then $\mut(G)=0$ if and only if $\mut(G_i)=0$ for at least one $i\in[k]$.
\end{corollary}

\section{Total mutual-visibility in Cartesian products}
\label{sec:cp}

In this section we consider the total mutual-visibility number of Cartesian product graphs. In the previous section we have seen that if $\mut(G) = 0$ or $\mut(H) = 0$, then $\mut(G\cp H) = 0$. Hence we may restrict our attention here to factor graphs with the total mutual-visibility number at least $1$.

To give general bounds we need the following concept. The \textit{independent total mutual-visibility number} $\muit(G)$ of $G$ is the cardinality of a largest independent total mutual-visibility set. Setting $\ell(G)$ to be the number of leaves of $G$, it follows from definitions that for any graph $G$, $\ell(G) \le \muit(G) \le \min \{\mut(G),\alpha(G)\}$. From Propositions~\ref{prop:for trees} and~\ref{prop:for block graphs} we know that the leaves set $L$ is a total mutual-visibility set of a tree $T$ and $\mut(G)=|L|$. Hence, if $n(T)\geq 3$, then $\mut(T)=|L|=\muit(T)$. Note in addition that $\muit(K_n) = 1$ while $\mut(K_n) = n$.

\begin{theorem}
\label{thm:cp-bounds}
If $G$ and $H$ are graphs of order at least $2$, $\muit(G) \ge 1$, and $\muit(H)\ge 1$, then
$$\max\{\muit(H)\mut(G), \muit(G)\mut(H)\} \le \mut(G\cp H) \le \min \{\mut(G) n(H), \mut(H) n(G)\}\,.$$
\end{theorem}

\proof
Let $I_G$ be an independent total mutual-visibility set of $G$ with $|I_G| = \muit(G)$,   and let $X_H$ be a $\mut$-set of $H$. Set $U = I_G \times X_H$. We claim that $U$ is a total mutual-visibility set of $G\cp H$.

Let $(g, h)$ and $(g', h')$ be arbitrary, different vertices of $G\cp H$. If $(g,h)$ and $(g',h')$ lie in the same $G$-layer or the same $H$-layer, then $x$ and $y$ are $U$-visible as layers are convex subgraphs of the product. Hence assume in the rest that $g\ne g'$ and $h\ne h'$.  

Let $g=g_0,g_1,\ldots,g_k=g'$ be the consecutive vertices of a shortest $g,g'$-path in $G$ whose internal vertices are not in $I_G$. Similarly, let $h=h_0,h_1,\ldots,h_\ell = h'$ be the consecutive vertices of a shortest $h,h'$-path in $H$ whose internal vertices are not in $X_H$. Assume first that $k=1$, that is, $gg'\in E(G)$. If $g\notin I_G$, then the path $(g,h)=(g_0,h_0),(g_0,h_1),\ldots,(g_0,h_\ell), (g_1,h_\ell)=(g',h')$ demonstrates that $(g, h)$ and $(g', h')$ are $U$-visible. If $g\in I_G$, then $g'\notin I_G$ and then the path $(g,h)=(g_0,h_0),(g_1,h_0),(g_1,h_1),\ldots,(g_1,h_\ell)=(g',h')$ demonstrates that $(g, h)$ and $(g', h')$ are $U$-visible. Assume in the following that $k\ge 2$. Consider now the path $P$ with the consecutive vertices
$$(g,h)=(g_0,h_0),(g_1,h_0),(g_1,h_1),\ldots,(g_1,h_\ell), (g_2,h_\ell),\ldots,(g_k,h_\ell)=(g',h')\,.$$
The path $P$ is a shortest $(g,h),(g',h')$-path with no internal vertex in $U$. Note that in all the above cases it is possible that $\ell = 1$ which happens if $hh'\in E(H)$. 

This proves the claim which in turn implies that $\mut(G\cp H) \ge |U| = \muit(G) \mut(H)$. By the commutativity of the Cartesian product, $\mut(G\cp H) \ge  \muit(H) \mut(G)$ and the lower bound follows.

Let $X$ be a total mutual of $G\cp H$. Since each $G$-layer $G^h$ is convex in $G\cp H$ we have $|X\cap V(G^h)| \le \mut(G)$, hence $\mut(G\cp H) \le \mut(G) n(H)$. Analogously $\mut(G\cp H) \le \mut(H) n(G)$.
\qed

In the rest of the section we give several exact results on $\mut(G\cp H)$ which also demonstrate that the bounds of Theorem~\ref{thm:cp-bounds} can be attained. We begin with the following sharpness result for the lower bound.

\begin{proposition}
\label{prop:both-factors-mut-1}
If $\mut(G\cp H) = 1$,  then $\mut(G) = 1$ and $\mut(H) = 1$.
\end{proposition}

\proof
Suppose first that $\mut(G)\ge 2$. Let $X_G$ be a total mutual-visibility set of $G$ with $|X_G|\ge 2$ and let $g_1, g_2\in X_G$. Let $X_H$ be a total mutual-visibility set of $H$ with $|X_H|\ge 1$ and let $h\in X_H$. Then we claim that $X = \{(g_1,h), (g_2,h)\}$ is a total mutual-visibility set of $G\cp H$. Indeed, if two vertices of $G\cp H$ are not lie in the layers that contain $(g_1,h)$ and $(g_2,h)$, then they are clearly $X$-visible. While if at least one of two vertices of $G\cp H$ lie in some of the layers in which $(g_1,h)$ and $(g_2,h)$, then the two vertices are again $X$-visible having Proposition~\ref{prop:subsets-are-ok} in mind. We conclude that $\mut(G) = 1$ and $\mut(H) = 1$.
\qed

\begin{figure}[ht!]
\begin{center}
\begin{tikzpicture}[scale=1.0,style=thick]
\tikzstyle{every node}=[draw=none,fill=none]
\def\vr{3pt} 

\begin{scope}[yshift = 0cm, xshift = 0cm]
 \node [below=0.5mm] at (0,0) {$x_1$};
    \node [below=0.5mm] at (1.5,0) {$x_2$};
    \node [below=0.5mm] at (3,0) {$x_3$};
    \node [above=0.5mm] at (0,1.5) {$x_4$};
    \node [below=0.5mm] at (1.5,1.5) {$x_5$};
    \node [above=0.5mm] at (3,1.5) {$x_6$};
    \node [above=0.5mm] at (1.5,2.5) {$x_7$};

\path (0,0) coordinate (x1);
\path (1.5,0) coordinate (x2);
\path (3,0) coordinate (x3);
\path (0,1.5) coordinate (x4);
\path (1.5,1.5) coordinate (x5);
\path (3,1.5) coordinate (x6);
\path (1.5,2.5) coordinate (x7);

\draw (x1) -- (x2) -- (x3);
\draw (x1) -- (x4) -- (x5)-- (x6) -- (x3);
\draw (x4) -- (x7) -- (x6);

\draw (x1)  [fill=white] circle (\vr);
\draw (x2)  [fill=white] circle (\vr);
\draw (x3)  [fill=white] circle (\vr);
\draw (x4)  [fill=white] circle (\vr);
\draw (x5)  [fill=white] circle (\vr);
\draw (x6)  [fill=white] circle (\vr);
\draw (x7)  [fill=white] circle (\vr);
\end{scope}
\end{tikzpicture}
\end{center}
\caption{The theta graph $\Theta(2,2,4)$. }
\label{fig:theta-example}
\end{figure}
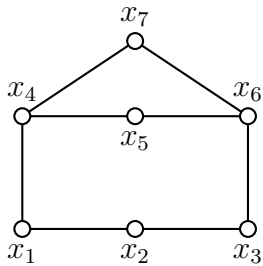

The converse of Proposition~\ref{prop:both-factors-mut-1} does not hold. For instance, consider the theta graph $\Theta(2,2,4)$ as presented in Fig.~\ref{fig:theta-example}.

It is straightforward to see that $\mut(\Theta(2,2,4)) = 1$ and that $\{x_5\}$ and $\{x_7\}$ are $\mut$-sets. In the Cartesian product $\Theta(2,2,4)\cp \Theta(2,2,4)$ one can see that $\{(x_5,x_5),(x_7,x_7)\}$ is a total mutual-visibility set, hence $\mut(\Theta(2,2,4)\cp \Theta(2,2,4))\ge 2$.

\cite[Corollary 3.7]{cicerone-2023} asserts that $\mu(K_n\cp K_m) = z(n,m;2,2)$, where $z(n,m;2,2)$ is the Zarankiewitz's number. To determine $z(n,m;2,2)$ is a notorious open problem~\cite{west-2021, z-1951}. Interestingly, the total mutual-visibility number of Cartesian products of complete graphs can be determined as follows, which further demonstrates that the lower bound of Theorem~\ref{thm:cp-bounds} is sharp.

\begin{proposition}
\label{prop:cp-complete-by-complete}
If $n, m\ge 2$, then $\mut(K_n\cp K_m) = \max \{n,m\}$.
\end{proposition}

\proof
Note first that the vertices of a single $K_n$-layer (or a single $K_m$-layer) form a total mutual-visibility set. Hence $\mut(K_n\cp K_m) \ge \max \{n,m\}$. To prove the other inequality, set $V(K_n) = [n]$ and $V(K_m) = [m]$, so that $V(K_n\cp K_m) = [n] \times [m]$.
Suppose without loss of generality that $n\ge m$. Let $U$ be an arbitrary total mutual-visibility set of $G\cp H$. If each $K_n$-layer contains at most one vertex of $U$ there is nothing to be proved. Assume now that some $K_n$-layer contains (at least) two vertices of $U$. By the symmetry we may assume that $(1,1)\in U$ and $(2,1)\in U$. We claim first that $(i,j)\notin U$, where $i, j\ge 2$. Indeed, if $(i,j)\in U$, then the vertices $(1,j)$ and $(i,1)$ are not $U$-visible. We claim second that $(1,j)\notin U$ for $2\le j\le m$. Indeed, if $(1,j)\in U$ for some $j\ge 2$, then the vertices $(1,1)$ and $(2,j)$ are not visible. We conclude that if $(1,1), (2,1)\in U$, then $U \subseteq V(K_n^1)$ and consequently $|U| \le n$. Analogously we see that if some $K_m$-layer contains (at least) two vertices of $U$, then $|U|\le m$. In any case, $\mut(K_n\cp K_m) \leq \max \{n,m\}$.
\qed

The next result (when $s\in \{3, 4\}$) also demonstrates sharpness of the lower bound of Theorem~\ref{thm:cp-bounds}.

\begin{proposition}
\label{thm:cp-cycle-by-complete graphs}
If $s\geq 3$ and $n\geq 3$, then
 $$\mut(C_s\cp K_n) = \left\{
\begin{array}{ll}
0; & s\geq 5,\\
n; & \mbox{otherwise}.\\
\end{array}\right.
$$
\end{proposition}
\proof
If $s\geq 5$, then Corollary~\ref{cor:girth} and Theorem~\ref{thm:cp} yield $\mut(C_s\cp K_n)=0$. In addition, $\mut(C_3\cp K_n) = n$ by Proposition~\ref{prop:cp-complete-by-complete}. Hence assume that $s = 4$ in the remaining proof.

By Theorem~\ref{thm:cp-bounds} we have $\mut(C_4\cp K_n) \ge n$. It remains to show that $\mut(C_4\cp K_n)\leq n$. Let $R$ be a $\mut$-set of $C_4\cp K_n$. If each $C_4$-layer contains at most one vertex of $R$, then $\mut(C_4\cp K_n)\leq n$ holds clearly. Suppose next that $R$ contains at least two vertices from the some $C_4$-layer. We may without loss of generality assume that $(1,n), (2,n)\in R$. Suppose that there exists another vertex $(i,j)\in R$, where $i\in[4]\setminus[2]$ and $j\in[n-1]$. If $i=3$, then the two vertices $(2,j)$ and $(3,n)$ are not $R$-visible. Similarly, the vertices $(1,j)$ and $(4,n)$ are not $R$-visible. This would thus mean that $|R| = 2$. We conclude that $\mut(C_4\cp K_n)=n$.
\qed

\begin{theorem}
\label{thm:cp-tree-by-graphs}
If $T$ is tree with $n(T)\ge 3$ and $H$ is a graph with $n(H)\ge 2$, then $\mut(T\cp H) = \mut(T) \mut(H)$.
\end{theorem}

\proof
The lower bound $\mut(T\cp H) \ge \mut(T) \mut(H)$ follows by Theorem~\ref{thm:cp-bounds} and the fact that $\mut(T) = \muit(T)$.

To prove that $\mut(T\cp H) \le \mut(T) \mut(H)$, consider an arbitrary $\mut$-set $R$ of $T\cp H$. Let $t\in V(T)$ be a vertex with $\deg_T(t) \ge 2$. Then $t$ is a non-bypass vertex of $T$. Hence, if $h\in V(H)$, then $(t,h)$ is a non-bypass vertex of $T\cp H$, thus $(t,h)\notin R$ by Lemma~\ref{lem:non-bypass-in-not-in}. Therefore, $R \cap V(^tH) = \emptyset$. So $R$ contains only vertices in $H$-layers corresponding to the leaves of $T$. Since each such layer can contain at most $\mut(H)$ vertices of $R$ we conclude that $\mut(T\cp H) \le |R| \leq \mut(T) \mut(H)$.
\qed

As a consequence of Theorem~\ref{thm:cp-tree-by-graphs} we obtain the following result which demonstrates sharpness of the upper bound of Theorem~\ref{thm:cp-bounds}.

\begin{corollary}
\label{cor:tree-by-complete}
If $T$ is tree with $n(T)\ge 3$, then $\mut(T\cp K_n) = n\cdot \mut(T)$.
\end{corollary}

\section{On the inequality $\mut(G\cp H) \le \mut(G) \mut(H)$}
\label{sec:cpII}

All the exact results obtained in Section~\ref{sec:cp} fulfil the bound
\begin{equation}
\label{eq:mut-times-mut}
\mut(G\cp H) \le \mut(G) \mut(H)\,.
\end{equation}
Hence one may wonder whether the upper bound of Theorem~\ref{thm:cp-bounds} can be improved/replaced by~\eqref{eq:mut-times-mut}. Before we answer the question, we prove another result where~\eqref{eq:mut-times-mut} holds.

A graph $G$ is a {\em generalized complete graph} if it is obtained by the join of an isolated vertex with a disjoint union of $k\ge 1$ complete graphs~\cite{tian-2021a}. We further say that $G$ is a {\em non-trivial generalized complete graph} if $k\ge 2$. Note that if $G$ is a non-trivial generalized complete graph, then $\mut(G) = n(G) - 1$.

\begin{theorem}
\label{thm:cp-generalized complete graphs}
If $G$ and $H$ are two non-trivial generalized complete graphs, then
$$\mut(G\cp H) \le \mut(G) \mut(H) = (n(G)-1)(n(H)-1)\,.$$
Moreover, the equality holds if and only if $G$ or $H$ is isomorphic to a star.
\end{theorem}

\proof
Let $V(G)=\{g_1,\ldots,g_{n(G)}\}$ and $V(H)=\{h_1,\ldots,h_{n(H)}\}$. Let $g_1$ and $h_1$ be the universal vertices of $G$ and $H$, respectively.

Note that $g_1$ is a non-bypass vertex of $G$ and $h_1$ is a non-bypass vertex of $H$. It follows that each vertex from the layer $^{g_1}H$ is a non-bypass vertex of $G\cp H$ as well as is each vertex from the layer $G^{h_1}$. Hence $G\cp H$ contains at least $n(G)+n(H)-1$ non-bypass vertices and so by~\eqref{eq:bounded-by-bp},
\begin{align*}
\mut(G\cp H) & \le \bp(G\cp H) \le n(G) n(H)- (n(G)+n(H)-1) \\
& =(n(G)-1)(n(H)-1)\,.
\end{align*}

To prove the equality case, by Theorem~\ref{thm:cp-tree-by-graphs} we know that $\mut(G\cp H) = (n(G)-1)(n(H)-1)$ if $G$ or $H$ is a star. Suppose in the rest that neither $G$ nor $H$ is a star. Then each of them contains an induced subgraph $K_3$. Without loss of generality, assume that $g_1, g_2, g_3$ induce a $K_3$ of $G$, and that $h_1, h_2, h_3$ induce a $K_3$ of $H$. Then the vertices $(g_2,h_2),(g_2,h_3),(g_3,h_3)$, and $(g_3,h_2)$ induce a $C_4$ of $G\cp H$. Since at most two vertices of this $C_4$ can lie in a total mutual-visibility set of $G\cp H$, we conclude that $\mut(G\cp H) < \bp(G\cp H)$.
\qed

In the rest we demonstrate that~\eqref{eq:mut-times-mut} does not hold in general. For this sake we say that a graph $G$ is {\em bypass over-visible} if it contains an independent bypass set of vertices $U$ which contains a $\mut$-set $U'$ as a {\bf proper} subset. Note that since $U'$ is an independent set, a bypass over-visible graph is a $(\muit,\mut)$-graph.

\begin{theorem}
\label{the:over-visible}
If $G$ and $H$ are bypass over-visible graphs, then $$\mut(G\cp H)>\mut(G) \mut(H).$$
\end{theorem}

\proof
Let $I_G$ and $I_H$ be largest independent bypass vertex sets of $G$ and $H$, which respectively contain $\mut$-sets $S_G$ and $S_H$ as proper subsets. Hence there exist vertices $u\in I_G \setminus S_G$ and $v\in I_H\setminus S_H$. We set $U = (S_G\times S_H) \cup \{(u,v)\}$ and claim that $U$ is a total mutual-visibility set of $G\cp H$.

Consider two arbitrary vertices $x=(g,h)$ and $y=(g',h')$ from $G\cp H$. Suppose first that $g=g'$. If $g\in S_G$, then $x$ and $y$ are $U$-visible because $S_H$ is a mutual visibility set of $H$. If $g=u$, then $x$ and $y$ are $U$-visible by Lemma~\ref{lem:bypass-vertex-is-good}. In all the other cases $V(^g{H}) \cap U = \emptyset$, hence there is nothing to prove. If $h=h'$, then $x$ and $y$ are $U$-visible by the same argument.

Assume in the rest that $g\ne g'$ and $h\neq h'$. Let $P_G: g=g_0,g_1,\ldots,g_k=g'$ be a shortest $g,g'$-path in $G$ whose internal vertices are not in $S_G$.
Similar, let $P_H: h=h_0,h_1,\ldots,h_\ell = h'$ be a shortest $h,h'$-path in $H$ whose internal vertices are not in $S_H$. The copy of $P_G$ in the layer $G^w$ will be denoted by $P_G^w$ and the copy of $P_H$ in the layer $^zH$ will be denoted by $^zP_H$.

Consider first the case $k=1$, that is, when $gg'\in E(G)$. If $(g_1,h_0)\notin U$, then the concatenation of the edge $(g_0,h_0)(g_1,h_0)$ and the path $^{g_1}P_H$ is a required $x,y$-path. Suppose next that $(g_1,h_0)\in U$. Then $h_0\in I_H$ and since $I_H$ is independent, we infer that $(g_0,h_\ell)\notin U$. Then the path $^{g_0}P_H$ followed by the edge $(g_0,h_\ell)(g_1,h_\ell)$ is a path which ensures that $x$ and $y$ are $U$-visible. Similarly we see that $x$ and $y$ are $U$-visible if $\ell = 1$. Note that the argument also applies when $ k = \ell = 1$.

We are left with the case when $k\ge 2$ and $\ell \ge 2$. Assume first that $u\neq g_i$ for $i\in[k-1]$. Then the vertices
$$(g,h)=(g_0,h_0),(g_1,h_0),(g_1,h_1),\ldots,(g_1,h_{\ell}),(g_2,h_{\ell}),\ldots,(g_k,h_\ell)=(g',h')$$
induce a shortest $g,g'$-path and the internal vertices of it are not in $U$. Hence $x$ and $y$ are $U$-visible. Similarly we see that $x$ and $y$ are $U$-visible if $v\neq h_j$ for $j\in[\ell-1]$. The remaining case is that $u=g_i$ for some $i\in[k-1]$ and $v = h_j$ for some $j\in[\ell-1]$. Then the vertices
$$(g,h)=(g_0,h_0),\ldots,(g_{i-1},h_0),(g_{i-1},h_1),\ldots,(g_{i-1},h_{\ell}),(g_i,h_{\ell}),\ldots,(g_k,h_\ell)=(g',h')$$
induce a shortest $x,y$-path in $G\cp H$ with no internal vertices in $U$. We conclude that in any case $x$ and $y$ are $U$-visible.
\qed

For $i\ge 2$ and $k\ge 3$, let $\Theta_i$ denote the theta graph $\Theta(p_1,\ldots, p_{k})$, where $$2 = p_1=\cdots = p_i < p_{i+1} \le \cdots \le p_k\,.$$ 
Let $a$ and $b$ be the vertices of $\Theta_i$ of degree $k$. Then $\BP(\Theta_i)$ consists of the degree $2$ vertices which are adjacent to $a$ and to $b$, so that $\bp(\Theta_i) = i$. Note in addition that $\BP(\Theta_i)$ is an independent set. On the other hand, $\BP(\Theta_i)$ is not a total mutual-visibility set, but becomes such a set if an arbitrary vertex is removed from it. Hence $\mut(\Theta_i) = i - 1$. We conclude that $\Theta_i$ is a bypass over-visible graph, hence any two such graphs fulfil the assumptions of Theorem~\ref{the:over-visible}.

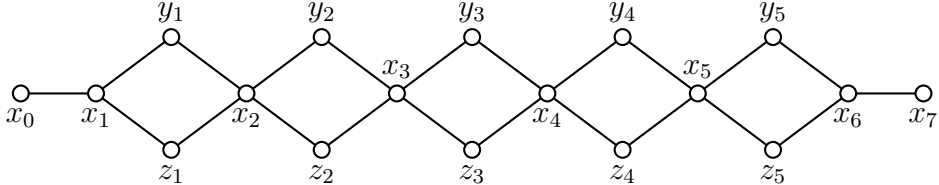
\begin{figure}[ht!]
\begin{center}
\begin{tikzpicture}[scale=1.0,style=thick]
\tikzstyle{every node}=[draw=none,fill=none]
\def\vr{3pt} 

\begin{scope}[yshift = 0cm, xshift = 0cm]
 \node [below=0.5mm] at (0,0) {$x_0$};
    \node [below=0.5mm] at (1,0) {$x_1$};
    \node [below=0.5mm] at (3,0) {$x_2$};
    \node [above=0.5mm] at (5,0) {$x_3$};
    \node [below=0.5mm] at (7,0) {$x_4$};
    \node [above=0.5mm] at (9,0) {$x_5$};
    \node [below=0.5mm] at (11,0) {$x_6$};
    \node [below=0.5mm] at (12,0) {$x_7$};
    \node [above=0.5mm] at (2,0.75) {$y_1$};
    \node [above=0.5mm] at (4,0.75) {$y_2$};
    \node [above=0.5mm] at (6,0.75) {$y_3$};
    \node [above=0.5mm] at (8,0.75) {$y_4$};
    \node [above=0.5mm] at (10,0.75) {$y_5$};
     \node [below=0.5mm] at (2,-0.75) {$z_1$};
    \node [below=0.5mm] at (4,-0.75) {$z_2$};
    \node [below=0.5mm] at (6,-0.75) {$z_3$};
    \node [below=0.5mm] at (8,-0.75) {$z_4$};
    \node [below=0.5mm] at (10,-0.75) {$z_5$};

\path (0,0) coordinate (x1);
\path (1,0) coordinate (x2);
\path (3,0) coordinate (x3);
\path (5,0) coordinate (x4);
\path (7,0) coordinate (x5);
\path (9,0) coordinate (x6);
\path (11,0) coordinate (x7);
\path (12,0) coordinate (x8);
\path (2,0.75) coordinate (x9);
\path (4,0.75) coordinate (x10);
\path (6,0.75) coordinate (x11);
\path (8,0.75) coordinate (x12);
\path (10,0.75) coordinate (x13);
\path (2,-0.75) coordinate (x14);
\path (4,-0.75) coordinate (x15);
\path (6,-0.75) coordinate (x16);
\path (8,-0.75) coordinate (x17);
\path (10,-0.75) coordinate (x18);
\draw (x1)-- (x2)--(x9)--(x3) -- (x10)--(x4)--(x11)--(x5)-- (x12)--(x6)--(x13)--(x7)--(x8);
\draw (x2)--(x14)-- (x3) -- (x15)--(x4)-- (x16) -- (x5)--(x17)-- (x6) -- (x18)--(x7);

\draw (x1)  [fill=white] circle (\vr);
\draw (x2)  [fill=white] circle (\vr);
\draw (x3)  [fill=white] circle (\vr);
\draw (x4)  [fill=white] circle (\vr);
\draw (x5)  [fill=white] circle (\vr);
\draw (x6)  [fill=white] circle (\vr);
\draw (x7)  [fill=white] circle (\vr);
\draw (x8)  [fill=white] circle (\vr);
\draw (x9)  [fill=white] circle (\vr);
\draw (x10)  [fill=white] circle (\vr);
\draw (x11)  [fill=white] circle (\vr);
\draw (x12)  [fill=white] circle (\vr);
\draw (x13)  [fill=white] circle (\vr);
\draw (x14)  [fill=white] circle (\vr);
\draw (x15)  [fill=white] circle (\vr);
\draw (x16)  [fill=white] circle (\vr);
\draw (x17)  [fill=white] circle (\vr);
\draw (x18)  [fill=white] circle (\vr);

\end{scope}
\end{tikzpicture}
\end{center}
\caption{The graph $G_5$. }
\label{fig:the graphs}
\end{figure}

Let $m\ge 1$. Then we define the the graph $G_m$ as follows. The vertex set is 
$$V(G_m) = \{x_0,x_1,\ldots,x_{m+2}\} \cup \{y_1,\ldots,y_m\} \cup \{z_1,\ldots,z_m\}\,.$$ 
For any $i\in [m]$ we connect $y_i$ and $z_i$ with $x_i$ and $x_{i+1}$. Finally add the edges $x_0x_1$ and $x_{m+1}x_{m+2}$. See Fig.~\ref{fig:the graphs} for $G_5$. 

It is straightforward to see that $\BP(G_m) = \{x_0, x_{m+2}\} \cup \{y_1,\ldots,y_m\} \cup \{z_1,\ldots,z_m\}$, hence $\bp(G_m)=2m+2$ and $\BP(G_m)$ is an independent set. In addition, since for any $i$, the vertices $y_i$ and $z_i$ cannot both lie in a total mutual-visibility set, the set $\{x_0,x_{m+2}, y_1,\ldots,y_{m}\}$ is an independent $\mut$-set of $G$. Hence $G_m$ is a bypass over-visible graph. Moreover, $|\BP(G_m)| - \mut(G_m) = m$. Now, by a parallel construction as in the proof of Theorem~\ref{the:over-visible} we find out that $\mut(G_m\cp G_m) \ge (m+2)^2 + m$, hence $\mut(G\cp H)$ can be arbitrary larger than $\mut(G) \mut(H)$.

\section{Concluding remarks}
\label{sec:conclude}

There are several possibilities how to continue the investigation of this paper, here we emphasize some of them. 

We have characterized the graphs $G$ with $\mut(G) = 0$. The next step would be to characterize the graphs $G$ with $\mut(G) = 1$ (and maybe also with $\mut(G) = 2$). 

In view of~\eqref{eq:bounded-by-bp} it would be interesting to consider the graphs $G$ with $\mut(G) = \bp(G)$. 

In this paper we had a closer look to the total mutual-visibility number of Cartesian product graphs. Some other graph operations also appear interesting for such investigations, in particular the strong product and the lexicographic product.

\section*{Acknowledgements}

Sandi Klav\v{z}ar was supported by the Slovenian Research Agency (ARRS) under the grants P1-0297, J1-2452, and N1-0285.
Jing Tian was supported by the Postgraduate Research Practice Innovation Program of Jiangsu Province No.KYCX22\_0323 and the Interdisciplinary Innovation Fund for Doctoral Students of Nanjing University of Aeronautics and Astronautics No.KXKCXJJ202204.


\end{document}